\documentclass[12pt]{article}
\usepackage{amssymb}

\usepackage{amsmath}

\newtheorem{theorem}{Theorem}

\newtheorem{definition}[theorem]{Definition}

\newtheorem{lemma}[theorem]{Lemma}

\newtheorem{proposition}[theorem]{Proposition}

\newenvironment{proof}[1][Proof]{\textbf{#1.} }{\ \rule{0.5em}{0.5em}}

\def \eref#1{\hbox{(\ref{#1})}}

\begin{document}

\title{Some Processes Associated with Fractional Bessel Processes}
\author{Yaozhong Hu\thanks{%
Y. Hu is supported in part by the National Science Foundation under Grant No.
DMS0204613} \\
Department of Mathematics\thinspace ,\ University of Kansas\\
405 Snow Hall\thinspace ,\ Lawrence, Kansas 66045-2142\\
Email: hu@math.ku.edu \and David Nualart\thanks{%
D. Nualart is supported in part by the MCyT Grant BFM2000-0598} \\
Facultat de Matem\`{a}tiques\thinspace , \ Universitat de Barcelona\\
Gran Via 585\thinspace\ \ 08007 Barcelona, Spain\\
Email: nualart@mat.ub.es}
\date{}
\maketitle

\begin{abstract}
Let $B=\{(B_{t}^{1},\ldots ,B_{t}^{d})\,,t\geq 0\}$ be a $d$-dimensional
fractional Brownian motion with Hurst parameter $H$ and let $R_{t}=%
\sqrt{(B_{t}^{1})^{2}+\cdots +(B_{t}^{d})^{2}}$ be the fractional Bessel
process. It\^{o}'s formula for the fractional Brownian motion leads to the
equation
$
R_{t}=\sum_{i=1}^{d}\int_{0}^{t}\frac{B_{s}^{i}}{R_{s}}%
dB_{s}^{i}+H(d-1)\int_{0}^{t}\frac{s^{2H-1}}{R_{s}}ds\,.
$
In the Brownian
motion case ($H=1/2$), $X_{t}=\sum_{i=1}^{d}\int_{0}^{t}\frac{B_{s}^{i}}{%
R_{s}}dB_{s}^{i}$ is a Brownian motion. In this paper it is shown that
$X_{t}$ is \underbar{not} a fractional Brownian motion if $H\not=1/2$. We
will study some other properties of this stochastic process as well.


\end{abstract}

\section{Introduction}

Let $B=\{(B_{t}^{1},\ldots ,B_{t}^{d})\,,t\geq 0\}$ be a $d$-dimensional
fractional Brownian motion with Hurst parameter $H\in (0,1)$. That is, the
components of $B$ are independent one-dimensional fractional Brownian
motions with Hurst parameter $H\in (0,1)$.

Denote the fractional Bessel process by $R_{t}=\sqrt{(B_{t}^{1})^{2}+\cdots+ 
(B_{t}^{d})^{2}}$. In the standard Brownian motion case there is an  extensive
literature on this process, see for example \cite{PY}. It is natural and
interesting to study this process for any other parameter $H$. If $d\geq 2$
and $1/2<H<1$, using
the It\^{o}'s formula for the fractional Brownian motion we obtain
\begin{equation}
R_{t}=\sum_{i=1}^{d}\int_{0}^{t}\frac{B_{s}^{i}}{R_{s}}dB_{s}^{i}+H(d-1)%
\int_{0}^{t}\frac{s^{2H-1}}{R_{s}}ds  \label{e8}
\end{equation}%
and for $d=1$ we have
\begin{equation}
|B_{t}|=\int_{0}^{t}\mathrm{sign}(B_{t})dB_{t}+H\int_{0}^{t}\delta
_{0}(B_{s})s^{2H-1}ds,  \label{e9}
\end{equation}%
where $\delta _{0}$ is the Dirac delta function, and the stochastic
integrals are interpreted in the divergence form. Equation (\ref{e8}) have
been proved in \cite{GN} in the case $H>\frac{1}{2}$, and for Equation (\ref%
{e9}) we refer to \cite{ChN}, \cite{GN}, \cite{HO} and \cite{HOS}.

In the classical Brownian motion case it is well-known from the L\'{e}vy's
characterization theorem that the first term in  the decomposition \eref{e8}
\begin{equation*}
X_{t}=%
\begin{cases}
\sum_{i=1}^{d}\int_{0}^{t}\frac{B_{s}^{i}}{R_{s}}dB_{s}^{i} & \qquad \mathrm{%
when}\ \ d\geq 2\cr\cr\int_{0}^{t}\mathrm{sign}(B_{t})dB_{t} & \qquad
\mathrm{when}\ \ d=1\cr%
\end{cases}%
\end{equation*}%
is a classical Brownian motion. It is then natural and interesting to ask
whether for any other $H$, the process $X=$ $\{X_{t}$, $t\geq 0\}$ is a fractional
Brownian motion or not.  The difficulty is that there is no characterization
as convenient as L\'{e}vy's one for general fractional Brownian motion ($%
H\not=1/2$). It is then difficult to show whether a stochastic process is a
fractional Brownian motion or not. In this paper, we shall prove that if $%
H\not=1/2$, then $\{X_{t}$, $t\geq 0\}$  is NOT a fractional Brownian
motion. Our approach to show this fact is based on the Wiener chaos
expansion (see for example \cite{DHP} and \cite{HO}). It seems to be the
natural method to be used here since there is no other powerful tool
available.

Although $\{X_{t}$, $t\geq 0\}$  is not a fractional Brownian motion, it
enjoys some properties that the fractional Brownian motion has, such as
self-similarity and long range dependence ($H>2/3$).   We will study these
and some other properties of the process $X$.

Section \ 2 will recall some preliminary results. Section 3 will study the
case $d=1$, namely,  the process $\int_{0}^{t}\mathrm{sign}(B_{t})dB_{t}$ and
Section 4 is devoted to the study of general dimension, ie, the process $%
\sum_{i=1}^{d}\int_{0}^{t}\frac{B_{s}^{i}}{R_{s}}dB_{s}^{i}$.

\section{Preliminaries}

Let $B=\{B_{t},t\geq 0\}$ be a fractional Brownian motion with Hurst
parameter $H\in (0,1)$. That is, $B$ is a zero mean Gaussian process with
the covariance function%
\begin{equation*}
\ R_{H}\left( t,s\right) =E(B_{t}B_{s})=\frac{1}{2}\left(
s^{2H}+t^{2H}-|t-s|^{2H}\right) .
\end{equation*}%
We denote by $K_{H}(t,s)$ the square integrable kernel such that%
\begin{equation*}
R_{H}(t,s)=\int_{0}^{t\wedge s}K_{H}(t,u)K_{H}(s,u)du.
\end{equation*}

Fix a time interval $[0,T]$, and let $\mathcal{H}$ be Hilbert space defined
as the closure of the set of step functions on $[0,T]$ with respect to the
scalar product
\begin{equation*}
\left\langle \mathbf{1}_{\left[ 0,t\right] },\mathbf{1}_{\left[ 0,s\right]
}\right\rangle _{\mathcal{H}}=R_{H}(t,s).\
\end{equation*}%
The mapping $\mathbf{1}_{\left[ 0,t\right] }\longrightarrow B_{t}$ can be
extended to an isometry between $\mathcal{H}$ and the Gaussian space $%
H_{1}(B)$ associated with $B$. We will denote this isometry by $\varphi
\longrightarrow B(\varphi )$.

The operator defined by

\begin{equation*}
\left( K_{H}^{\ast }\mathbf{1}_{[0,t]}\right) (s)=K_{H}(t,s)\mathbf{1}%
_{[0,t]}(s).
\end{equation*}%
can be extended to a linear isometry between $\mathcal{H}$ and $L^{2}(0,T)$.
This operator can be expressed in terms of fractional operators. More
precisely, if $H>\frac{1}{2}$ we have%
\begin{equation*}
\left( K_{H}^{\ast }\varphi \right) (s)=c_{H}\Gamma (H-\frac{1}{2})s^{\frac{1%
}{2}-H}(I_{T-}^{H-\frac{1}{2}}u^{H-\frac{1}{2}}\varphi (u))(s)
\end{equation*}%
and if $H<\frac{1}{2}$
\begin{equation*}
\left( K_{H}^{\ast }\varphi \right) (s)=c_{H}\Gamma (H+\frac{1}{2})\ s^{%
\frac{1}{2}-H}(D_{T-}^{\frac{1}{2}-H}u^{H-\frac{1}{2}}\varphi (u))(s),
\end{equation*}%
where $c_{H}=\sqrt{\frac{2H}{(1-2H)\beta (1-2H,H+1/2)}}$, and for any $%
\alpha >0$ we denote by $I_{T-}^{\alpha }$  (resp. $D_{T-}^{\alpha }$) the
fractional integral (resp. derivative) operator given by%
\begin{equation*}
I_{T-}^{\alpha }f\left( t\right) =\frac{\left( -1\right) ^{-\alpha }}{\Gamma
\left( \alpha \right) }\int_{t}^{T}\left( s-t\right) ^{\alpha -1}f\left(
s\right) ds
\end{equation*}%
(resp.
\begin{equation*}
\left. D_{T-}^{\alpha }f\left( t\right) =\dfrac{\left( -1\right) ^{\alpha }}{%
\Gamma \left( 1-\alpha \right) }\left( \dfrac{f\left( t\right) }{\left(
T-t\right) ^{\alpha }}+\alpha \int_{t}^{T}\dfrac{f\left( t\right) -f\left(
s\right) }{\left( s-t\right) ^{\alpha +1}}ds\right) \right) .
\end{equation*}

We denote by $D$ and $\delta $ the derivative and divergence operators that
can be defined in the framework of the Malliavin calculus with respect to
the process $B$. Let $\mathbb{D}^{k,p}$, $p>1$, $k\in \mathbb{R}$, be the
corresponding Sobolev spaces. We recall that the divergence operator $\delta
$ is defined by means of the duality relationship
\begin{equation}
\mathbb{E}(F\delta (u))=E\left\langle DF,u\right\rangle _{\mathcal{H}},
\label{e5}
\end{equation}%
where $u$ is a random variable in $L^{2}(\Omega ;\mathcal{H})$. We say that $%
u$ belongs to the domain of the divergence, denoted by Dom $\delta $, if
there is a square integrable random variable $\delta (u)$ such that (\ref{e5}%
) holds for any $F\in $ $\mathbb{D}^{1,2}$.

The domain of the divergence operator is sometimes too small. For instance,
in \cite{ChN} it is proved that the process $u=B$ belongs to $L^{2}(\Omega ;%
\mathcal{H})$ if and only if $H>\frac{1}{4}$. On the other hand, in \cite%
{CNT} it is proved that for all $t\geq 0$, the process $\mathrm{sign}(B_{t})$
belongs to the domain of the divergence when $H>\frac{1}{3}$.

Following the approach of \cite{ChN} it is possible to extend the domain of
the divergence operator to processes whose trajectories are not necessarily
in the space $\mathcal{H}$. Set $\mathcal{H}_{2}=\left( K_{H}^{\ast }\right)
^{-1}\left( K_{H}^{\ast ,a}\right) ^{-1}(L^{2}(0,T))$, where $K_{H}^{\ast
,a} $denotes the adjoint of the operator $K_{H}^{\ast }$. Denote by $%
\mathcal{S}_{\mathcal{H}}$ the space of smooth and cylindrical random
variables of the form
\begin{equation}
F=f(B(\phi _{1}),\ldots ,B(\phi _{n})),  \label{f13}
\end{equation}%
where $n\geq 1$, $f\in C_{b}^{\infty }\left( \mathbb{R}^{n}\right) $ ($f$
and all its partial derivatives are bounded), and $\phi _{i}\in \mathcal{H}%
_{2}$. $\ $

\begin{definition}
\label{def1} Let $\ u=\{u_{t},t\in \lbrack 0,T]\}$ be a measurable process
such that%
\begin{equation*}
\mathbb{E}\left( \int_{0}^{T}u_{t}^{2}dt\right) <\infty .
\end{equation*}%
We say that $u\in \mathrm{Dom}_{T}^{\ast }\delta $ (extended domain of the
divergence in $[0,T]$) if there exists a random variable $\delta (u)\in
L^{2}(\Omega )$ such that for all $F\in \mathcal{S}_{\mathcal{H}}$ we have
\begin{equation*}
\int_{0}^{T}\mathbb{E}(u_{t}K_{H}^{\ast ,a}K_{H}^{\ast }D_{t}F)dt=\mathbb{E}%
(\delta (u)F).
\end{equation*}
\end{definition}

In \cite{ChN} it is proved that for any $H\in (0,1)$,  the process $\mathrm{%
sign}(B_{t})\mathbf{\ }$ belongs to the extended domain of the divergence in
any time interval $[0,T]$ and the following version of Tanaka's formula
holds {\small \ }%
\begin{equation}
|B_{t}|=\int_{0}^{t}\mathrm{sign}(B_{s})dB_{s}+H\int_{0}^{t}\delta
_{0}(B_{s})s^{2H-1}ds.  \label{e6}
\end{equation}%
(see also \cite{HO}, \cite{HOS} for this and a more general formula). In
this formula $L_{t}^{a}=H\int_{0}^{t}\delta _{0}(B_{s})s^{2H-1}ds$ is the
the density of the occupation measure%
\begin{equation*}
\Gamma \mapsto 2H\int_{0}^{t}1_{\Gamma }(B_{s})s^{2H-1}ds.
\end{equation*}

\section{The process $\protect\int_{0}^{t}\mathrm{sign}(B_{t})dB_{t}$}

Define the process $X=\{X_{t},t\geq 0\}$, by%
\begin{equation*}
X_{t}=\int_{0}^{t}\mathrm{sign}(B_{s})dB_{s}.
\end{equation*}

In the case of the classical Brownian motion ($H=\frac{1}{2}$), the process $%
X$ turns out to be a Brownian motion. We will show first that for any $H\in
(0,1)$,  $X$ is a $H$-self-similar process, that is, for all $a>0$ the
processes $\{X_{at},t\geq 0\}$ and  $\{a^{H}X_{t},t\geq 0\}$ have the same
law.

\begin{proposition}
The process $X=\{X_{t},t\geq 0\}$ is $H$-self-similar.
\end{proposition}

\begin{proof}
Using the self-similarity property of the fractional Brownian motion and
Tanaka's formula (\ref{e6}) yields that for any $a>0$
\begin{eqnarray*}
X_{at} &=&|B_{at}|-H\int_{0}^{at}\delta _{0}(B_{s})s^{2H-1}ds \\
&=&|B_{at}|-H\int_{0}^{t}\delta _{0}(B_{au})(au)^{2H-1}adu \\
&&\overset{d}{=}a^{H}|B_{t}|-a^{2H}H\int_{0}^{t}\delta
_{0}(a^{H}B_{u})u^{2H-1}du \\
&=&a^{H}X_{t},
\end{eqnarray*}%
where the symbol $\overset{d}{=}$ means that the distributions  of both
processes are  the same. This completes the proof.
\end{proof}

Then, it is natural to conjecture that for any $H$, the process $X_{t}$ is
a fractional Brownian motion of Hurst parameter $H$. We will see that this
is no longer true if $H\neq \frac{1}{2}$, although the process $X_{t}$
shares some of the properties of the fractional Brownian motion.

Let us first find the Wiener chaos expansion of the process $\mathrm{sign}%
(B_{t})$. We will denote by $I_{n}$ the multiple Wiener integral with
respect to the process $B$.

\begin{lemma}
\label{lem2} Let $0<H<1$. We have the following chaos expansion for $\mathrm{%
sign}(B_{t})$:
\begin{equation}
\mathrm{sign}(B_{t})=\sum_{\substack{ k=0}}^{\infty }b_{2k+1}I_{2k+1}(1),
\label{e.1.1}
\end{equation}%
where%
\begin{equation*}
b_{2k+1}=\frac{2(-1)^{k}}{(2k+1)\sqrt{2\pi }t^{(2k+1)H}k!2^{k}}.
\end{equation*}
\end{lemma}

\begin{proof}
Denote by $p_{\varepsilon }(x)=\frac{1}{\sqrt{2\pi \varepsilon }}%
e^{-x^{2}/\varepsilon }$, $x\in \mathbb{R}$, $\varepsilon >0$, the heat
kernel. The function%
\begin{equation*}
f_{\varepsilon }(x)=2\int_{-\infty }^{x}p_{\varepsilon }(y)dy-1
\end{equation*}%
converges to $\mathrm{sign}(x)$ as $\varepsilon $ tends to zero. Hence, $%
f_{\varepsilon }(B_{t})$ converges to $\mathrm{sign}(B_{t})$ in $%
L^{2}(\Omega )$ as $\varepsilon $ tends to zero. The application of
Stroock's formula yields%
\begin{equation}
f_{\varepsilon }(B_{t})=\sum_{n=0}^{\infty }a_{n}^{\varepsilon
}(t)\int_{0<s_{1}<\cdots <s_{n}<t}dB_{s_{1}}\cdots dB_{s_{n}},  \label{e3}
\end{equation}%
where%
\begin{eqnarray*}
a_{n}^{\varepsilon }(t) &=&\mathbb{E}\left[ D^{n}(f_{\varepsilon }(B_{t}))%
\right] =2\mathbb{E}\left[ p_{\varepsilon }^{(n-1)}(B_{t})\right]  \\
&=&2(-1)^{n-1}\frac{\partial ^{n-1}}{\partial y^{n-1}}\mathbb{E}\left[
p_{\varepsilon }(B_{t}-y)\right] |_{y=0} \\
&=&2(-1)^{n-1}p_{\varepsilon +t^{2H}}^{(n-1)}(0).
\end{eqnarray*}%
Taking the limit of (\ref{e3}) in $L^{2}(\Omega )$ as $\varepsilon $ tends
to zero we obtain%
\begin{equation*}
\mathrm{sign}(B_{t})=\sum_{n=0}^{\infty }a_{n}(t)\int_{0<s_{1}<\cdots
<s_{n}<t}dB_{s_{1}}\cdots dB_{s_{n}},
\end{equation*}%
where $a_{n}(t)=\lim_{\varepsilon \downarrow 0}a_{n}^{\varepsilon
}(t)=2(-1)^{n-1}p_{t^{2H}}^{(n-1)}(0)$. As a consequence, $a_{n}(t)=0$ if $n$
is even and%
\begin{equation*}
a_{n}(t)=\frac{2(-1)^{k}(2k)!}{\sqrt{2\pi }t^{nH}k!2^{k}}
\end{equation*}%
if $n=2k+1$.
\end{proof}

Using Stirling's formula we obtain
\begin{eqnarray*}
\mathbb{E}\left[ I_{2k+1}(b_{2k+1})\right] ^{2} &=&\frac{4(2k+1)!t^{2H+1}}{%
(2k+1)2\pi t^{2H+1}\left( k!2^{k}\right) ^{2}} \\
&=&\frac{4(2k)!}{(2k+1)2\pi \left( k!2^{k}\right) ^{2}} \\
&\backsimeq &Ck^{-3/2},
\end{eqnarray*}%
and we have proved the following proposition.

\begin{proposition}
For any $0<H<1$, the random variable $\mathrm{sign}(B_{t})$ belongs to the
Sobolev space $\mathbb{D}^{\alpha ,2}$ for any $\alpha <\frac{1}{2}$.
\end{proposition}

Now it is easy to obtain the chaos expansion of $\int_0^t \mathrm{sign}%
(B_{s})dB_s$.

\begin{proposition}
For any $0<H<1$,
\begin{equation*}
\int_{0}^{t}\mathrm{sign}(B_{s})dB_{s}=\sum_{k=1}^{\infty
}c_{k}I_{2k}(h_{2k})\,,
\end{equation*}%
where
\begin{equation*}
c_{k}=\frac{(-1)^{k-1}}{\sqrt{2\pi }(2k-1)(k-1)!2^{k-2}}
\end{equation*}%
and
\begin{equation*}
h_{2k}(s_{1},\ldots ,s_{2k})=\left( s_{1}\vee s_{2}\vee \cdots \vee
s_{2k}\right) ^{-(2k-1)H }\,.
\end{equation*}
\end{proposition}

A consequence of this proposition is the following

\begin{proposition}
For any $0<H<1$ and $t>0$, the random variable $\int_0^t \mathrm{sign}%
(B_{s})dB_s$ belongs to the Sobolev space $\mathbb{D}^{\alpha ,2}$ for any $%
\alpha <1/2$.
\end{proposition}

\begin{proof}
It is easy to check that there is a constant $C>0$ such that
\begin{equation*}
{{\mathbb{E}}}\left[ I_{2k}(h_{2k})\right] ^{2}\leq C\frac{(2k)!}{(2k-1)^{2}%
\left[ (k-1)!\right] ^{2}2^{2k}}\,.
\end{equation*}%
Therefore
\begin{eqnarray*}
{{\mathbb{E}}}\left[ c_k I_{2k}(h_{2k})\right] ^{2}
&\leq &C\frac{(2k)!}{\left( k!2^{k}\right) ^{2}} \\
&\leq &Ck^{-3/2}\,.
\end{eqnarray*}%
This proves the proposition.
\end{proof}

The next proposition states that $\int_{0}^{t}\mathrm{sign}(B_{t})dB_{t}$ is
not a fractional Brownian motion.

\begin{proposition}
The process $X=\{X_{t},t\geq 0\}$ is not a fractional Brownian motion.
\end{proposition}

\begin{proof}
Suppose that $X$ is a fractional Brownian motion. Then it
is a fractional Brownian motion with Hurst parameter $H$
since it is self-similar with parameter $H$.
Then, the process%
\begin{equation*}
Y_{t}=\int_{0}^{t}\eta _{H}(t,r)dX_{r}
\end{equation*}%
must be a standard Brownian motion with respect to the filtration generated
by $X$, where%
\begin{equation*}
\eta _{H}(t,r)=\left( K_{H}^{\ast }\right) ^{-1}(\mathbf{1}_{[0,t]})(r).
\end{equation*}%
We claim that
\begin{equation}
Y_{t}=Z_{t},  \label{e4}
\end{equation}%
where%
\begin{equation}
Z_{t}=\int_{0}^{t}\eta _{H}(t,r)\mathrm{sign}(B_{r})dB_{r}.  \label{e.1.2}
\end{equation}%
In fact, set $t_{k}^{n}=\frac{tk}{n}$, $k=0,\ldots ,n$, and consider the
approximations%
\begin{equation*}
Y_{t}^{n}=\sum_{k=1}^{n}\ \eta _{H}(t,t_{k-1}^{n})\left(
X_{t_{k}^{n}}-X_{t_{k-1}^{n}}\right) .
\end{equation*}%
We know that $Y_{t}^{n}$ converges in $L^{2}(\Omega )$ to $Y_{t}$, because
the functions%
\begin{equation*}
\sum_{k=1}^{n}\eta _{H}(t,t_{k-1}^{n})\mathbf{1}_{[t_{k-1}^{n},t_{k}^{n})}(r)
\end{equation*}%
converge to $\eta _{H}(t,r)\mathbf{1}_{[0,t)}(r)$ in the norm of the Hilbert
space $\mathcal{H}$. On the other hand, by Definition \ref{def1}, for any
smooth and cylindrical random variable $F\in \mathcal{S}_{\mathcal{H}}$ we
have
\begin{eqnarray*}
\mathbb{E}(FY_{t}^{n}) &=&\mathbb{E}\left( \left\langle
D_{r}F,\sum_{k=1}^{n}\eta _{H}(t,t_{k-1}^{n})\mathbf{1}%
_{[t_{k-1}^{n},t_{k}^{n})}(r)\mathrm{sign}(B_{r})\right\rangle _{\mathcal{H}%
}\right) \\
&=&\mathbb{E}\left( \left\langle \Gamma _{H,T}^{\ast ,a}\Gamma _{H,T}^{\ast
}D_{r}F,\sum_{k=1}^{n}\eta _{H}(t,t_{k-1}^{n})\mathbf{1}%
_{[t_{k-1}^{n},t_{k}^{n})}(r)\mathrm{sign}(B_{r})\right\rangle
_{L^{2}(0,T)}\right) .
\end{eqnarray*}%
As before this converges to%
\begin{eqnarray*}
&&\mathbb{E}\left( \left\langle K^{\ast ,a}K^{\ast }D_{r}F,\eta _{H}(t,r)%
\mathbf{1}_{[0,t)}(r)\mathrm{sign}(B_{r})\right\rangle _{L^{2}(0,T)}\right)
\\
&=&\mathbb{E}(FZ_{t}),
\end{eqnarray*}%
as $n$ tends to infinity. So (\ref{e4}) holds.

We can write, using Lemma \ref{lem2}%
\begin{equation*}
Z_{t}=\sum_{\substack{ k=0}}^{\infty }b_{k}\int_{0}^{t}\eta
_{H}(t,r)r^{-(2k+1)H}I_{2k+1}(\mathbf{1}_{[0,r]}^{\otimes (2k+1)})dB_{r},
\end{equation*}%
where%
\begin{equation}
b_{k}=\frac{(-1)^{k}}{\sqrt{2\pi }(2k+1)k!2^{k-1}}.  \label{e7}
\end{equation}%
So,%
\begin{equation*}
Z_{t}=\sum_{\substack{ k=0}}^{\infty }b_{k}I_{2k+2,t}(f_{2k+2}),
\end{equation*}%
where%
\begin{eqnarray*}
&&f_{2k+2}(t,s_{1},\ldots ,s_{2k+2}) \\
&=&\mathrm{symm}\left( \eta _{H}(t,s_{2k+2})s_{2k+2}^{-(2k+1)H}\mathbf{1}%
_{[0,2k+2]}^{\ }(s_{1})\cdots \mathbf{1}_{[0,2k+2]}^{\ }(s_{2k+1})\right) \\
&=&\frac{1}{2k+1} \eta _{H}(t,s_{1}\vee \cdots \vee s_{2k+2})\left( s_{1}\vee \cdots \vee
s_{2k+2}\right) ^{-(2k+1)H},
\end{eqnarray*}%
and $I_{2k+2,t}(f)$ denotes $I_{2k+2}(f\mathbf{1}_{[0,t]}^{\ \otimes
(2k+2)}) $. We can transform these multiple stochastic integrals into
integrals with respect to a standard Brownian motion, using the operator $%
K_{H}^{\ast }$. In this way we obtain%
\begin{equation*}
I_{2k+2,t}(f_{2k+2})=I_{2k+2,t}^{W}(K_{H}^{\ast \otimes (2k+2)}f_{2k+2}),
\end{equation*}%
and the process%
\begin{equation*}
Z_{t}=\sum_{\substack{ k=0}}^{\infty }b_{k}I_{2k+2,t}^{W}(K_{H}^{\ast
\otimes (2k+2)}f_{2k+2})
\end{equation*}%
is a Brownian motion with respect to the filtration generated by $W$. Hence,
every component of the chaos expansion is a martingale with respect to the
filtration generated by $W$. In particular, this implies that the
coefficient of the second
chaos $K_{H}^{\ast \otimes 2}f_{2}(t,s_{1},s_{2})$ must not depend on $t$.

For $H>\frac{1}{2}$ we have%
\begin{eqnarray*}
K_{H}^{\ast \otimes 2}f_{2}(t,s_{1},s_{2}) &=&d_{H}^{2}\ \left(
s_{1}s_{2}\right) ^{\frac{1}{2}-H} \\
&&\times \left[ I_{t-}^{\left( H-\frac{1}{2}\right) \otimes 2}\left(
u_{1}u_{2}\right) ^{-\frac{1}{2}}\eta _{H}(t,u_{1}\vee u_{2})\right]
(s_{1},s_{2})
\end{eqnarray*}%
where $d_{H}=c_{H}\Gamma (H-\frac{1}{2})$. We have used the fact that%
\begin{equation*}
\left( I_{t-}^{\alpha }f\right) \mathbf{1}_{[0,t]}=I_{T-}^{\alpha }\left( f%
\mathbf{1}_{[0,t]}\right) .
\end{equation*}%
Then%
\begin{eqnarray*}
&&\left[ I_{t-}^{\left( H-\frac{1}{2}\right) \otimes 2}\left(
u_{1}u_{2}\right) ^{-\frac{1}{2}}\eta _{H}(t,u_{1}\vee u_{2})\right]
(s_{1},s_{2}) \\
&=&\frac{1}{\Gamma (H-\frac{1}{2})^{2}}\int_{s_{2}}^{t}\int_{s_{1}}^{t}%
\left( \left( u_{1}-s_{1}\right) \left( u_{2}-s_{2}\right) \right) ^{-\frac{1%
}{2}} \\
&&\times \eta _{H}(t,\left( u_{1}-s_{1}\right) \vee \left(
u_{2}-s_{2}\right) )du_{1}du_{2}.
\end{eqnarray*}%
Taking $t=\max (s_{1},s_{2})$, we would have $K_{H}^{\ast \otimes
2}f_{2}(t,s_{1},s_{2})=0$, because%
\begin{equation*}
\eta _{H}(t,r)\leq Ct^{H-\frac{1}{2}}r^{H-\frac{1}{2}}(t-r)^{\frac{1}{2}-H}.
\end{equation*}%
Hence, $\eta _{H}(t,u_{1}\vee u_{2})=0$, which leads to a contradiction.

Suppose now that $H<\frac{1}{2}$. In this case we have%
\begin{eqnarray*}
K_{H}^{\ast \otimes 2}f_{2}(t,s_{1},s_{2}) &=&e_{H}^{2}\ \left(
s_{1}s_{2}\right) ^{\frac{1}{2}-H} \\
&&\times \left[ D_{t-}^{\left( \frac{1}{2}-H\right) \otimes 2}\left(
u_{1}u_{2}\right) ^{-\frac{1}{2}}\eta _{H}(t,u_{1}\vee u_{2})\right]
(s_{1},s_{2}),
\end{eqnarray*}%
where $e_{H}=c_{H}\Gamma (H+\frac{1}{2})$, using again that $\left(
D_{t-}^{\alpha }f\right) \mathbf{1}_{[0,t]}=D_{T-}^{\alpha }\left( f\mathbf{1%
}_{[0,t]}\right) $. Notice that%
\begin{equation*}
\eta _{H}(t,r)=\frac{1}{e_{H}\Gamma (\frac{1}{2}-H)}r^{\frac{1}{2}%
-H}\int_{r}^{t}(y-r)^{-\frac{1}{2}-H}y^{H-\frac{1}{2}}dy.
\end{equation*}%
As a consequence, $\eta _{H}(t,r)$ behaves as $Cr^{\frac{1}{2}-H}(t-r)^{%
\frac{1}{2}-H}.$We have%
\begin{eqnarray*}
&&\left[ D_{t-}^{\left( \frac{1}{2}-H\right) \otimes 2}\left(
u_{1}u_{2}\right) ^{-\frac{1}{2}}\eta _{H}(t,u_{1}\vee u_{2})\right]
(s_{1},s_{2}) \\
&=&\dfrac{1}{\Gamma \left( H+\frac{1}{2}\right) ^{2}}D_{t-}^{\frac{1}{2}%
-H}\left( \dfrac{\left( s_{1}s_{2}\right) ^{-\frac{1}{2}}\eta
_{H}(t,s_{1}\vee s_{2})}{\left( t-s_{1}\right) ^{\frac{1}{2}-H}\left(
t-s_{2}\right) ^{\frac{1}{2}-H}}\right. \\
&&+\left( \frac{1}{2}-H\right) \int_{s_{1}}^{t}\dfrac{\left(
s_{1}s_{2}\right) ^{-\frac{1}{2}}\eta _{H}(t,s_{1}\vee s_{2})-\left(
ys_{2}\right) ^{-\frac{1}{2}}\eta _{H}(t,y\vee s_{2})}{\left( y-s_{1}\right)
^{\frac{3}{2}-H}\left( t-s_{2}\right) ^{\frac{1}{2}-H}}dy \\
&&+\left( \frac{1}{2}-H\right) \int_{s_{2}}^{t}\dfrac{\left(
s_{1}s_{2}\right) ^{-\frac{1}{2}}\eta _{H}(t,s_{1}\vee s_{2})-\left(
ys_{1}\right) ^{-\frac{1}{2}}\eta _{H}(t,y\vee s_{1})}{\left( y-s_{2}\right)
^{^{\frac{3}{2}-H}}\left( t-s_{1}\right) ^{\frac{1}{2}-H}}dy \\
&&+\left( \frac{1}{2}-H\right) ^{2}\int_{s_{1}}^{t}\int_{s_{2}}^{t}\left(
y-s_{1}\right) ^{H-\frac{3}{2}}\left( z-s_{2}\right) ^{H-\frac{3}{2}}\left(
s_{1}s_{2}\right) ^{-\frac{1}{2}}\left[ \eta _{H}(t,s_{1}\vee s_{2})\right.
\\
&&\left. \left. -\left( ys_{2}\right) ^{-\frac{1}{2}}\eta _{H}(t,y\vee
s_{2})-\left( zs_{1}\right) ^{-\frac{1}{2}}\eta _{H}(t,z\vee s_{1})+\left(
yz\right) ^{-\frac{1}{2}}\eta _{H}(t,y\vee z)\right] dydz\right) .
\end{eqnarray*}%
Taking again $t=\max (s_{1},s_{2})$, we would have $K_{H}^{\ast \otimes
2}f_{2}(t,s_{1},s_{2})=0$ which leads to a contradiction.
\end{proof}

Consider the covariance between two increments of the process $X$:%
\begin{equation*}
r(n):={{\mathbb{E}}}\left[ \left( X_{a+1}-X_{a}\right) (X_{n+1}-X_{n})\right]
,
\end{equation*}%
where $0<a\leq n$. We say that $X$ is long-range dependent if for any $a>0$,%
\begin{equation*}
\sum_{n\geq a}\left| r(n)\right| =\infty .
\end{equation*}

The next proposition studies the long-range dependence properties of the
process $X$. We see that this property differs from that of fractional
Brownian motion.

\begin{proposition}
Let $X_{t}=\int_{0}^{t}\mathrm{sign}(B_{s})dB_{s}$. If $H\geq 2/3$, then $%
X_{t}$ is long-range dependent and if $1/2<H<2/3$, then $X_{t}$ is not
long-range dependent.
\end{proposition}

\begin{proof}
From Lemma \ref{lem2} we can deduce the Wiener chaos expansion of the random
variable $X_{t}$. In fact, we have
\begin{equation*}
X_{t}=\sum_{k=0}^{\infty }b_{k}I_{2k+2,t}(h_{2k+2})\,,
\end{equation*}%
where $b_{k}$ is defined in (\ref{e7}) and
\begin{equation*}
h_{2k+2}(s_{1},\ldots ,s_{2k+2})=\left( s_{1}\vee \cdots \vee
s_{2k+2}\right) ^{-(2k+1)H}\,.
\end{equation*}%
Let us compute the covariance of $X_{s}$ and $X_{t}-X_{r}$, where $0<r<t$.
From the It\^{o} isometry of multiple stochastic integrals it follows that
\begin{equation*}
{{\mathbb{E}}}\left[ X_{s}(X_{t}-X_{r})\right] =\sum_{k=0}^{\infty }b_{k}^{2}%
{{\mathbb{E}}}\left[ I_{2k+2,s}(h_{2k+2})\left[
I_{2k+2,t}(h_{2k+2})-I_{2k+2,r}(h_{2k+2})\right] \right] .
\end{equation*}%
We have
\begin{eqnarray*}
{{\mathbb{E}}}\left[ X_{s}(X_{t}-X_{r})\right] &=&\sum_{k=0}^{\infty
}b_{k}^{2}(2k+2)!\left\langle h_{2k+2}\mathbf{1}_{[0,s]}^{\otimes
(k+2)},h_{2k+2}\left( \mathbf{1}_{[0,t]}^{\otimes (k+2)}-\mathbf{1}%
_{[0,r]}^{\otimes (k+2)}\right) \right\rangle _{\mathcal{H}^{2k+2}} \\
&\geq &2b_{0}^{2}\left\langle h_{2}\mathbf{1}_{[0,s]}^{\otimes
2},h_{2}\left( \mathbf{1}_{[0,t]}^{\otimes 2}-\mathbf{1}_{[0,r]}^{\otimes
2 }\right) \right\rangle _{\mathcal{H}^{2}} \\
&\geq &\frac{b_{0}^{2}}{2}\int_{r}^{t}\int_{0}^{r}\int_{{0}}^{s}\int_{{0}%
}^{s}s_{2}^{-H}t_{2}^{-H}\phi (s_{1},t_{1})\phi
(s_{2},t_{2})ds_{1}ds_{2}dt_{1}dt_{2},
\end{eqnarray*}%
where $\phi (s,t)=H(2H-1)|t-s|^{2H-2}$.

Thus let $s=1$, $r=n$, and $t=n+1$ and we have
\begin{eqnarray*}
r(n)&:= &{{\mathbb{E}}}\left[ \left( X_{a+1}-X_{a}\right) (X_{n+1}-X_{n})%
\right] \\
&\geq
&C\int_{n}^{n+1}\int_{0}^{n}\int_{a}^{a+1}%
\int_{a}^{a+1}s_{2}^{-H}t_{2}^{-H}(t_{1}-s_{1})^{2H-2}(t_{2}-s_{2})^{2H-2}ds_{1}ds_{2}dt_{1}dt_{2}
\\
&\geq
&C\int_{n}^{n+1}\int_{a+2}^{n}\int_{a}^{a+1}%
\int_{a}^{a+1}s_{2}^{-H}t_{2}^{-H}(t_{1}-a-1)^{2H-2}(t_{2}-a-1)^{2H-2}ds_{1}ds_{2}dt_{1}dt_{2}
\\
&\approx &Cn^{3H-3}\,,
\end{eqnarray*}%
as $n$ tends to infinity. Thus if $H\geq 2/3$, $\sum_{n\geq a}r(n)=\infty $.

If $H<2/3$, then we use another approach. Set, as before%
\begin{equation*}
r(n)=\mathbb{E}\left[ \left( \int_{a}^{a+1}\mathrm{sign}B_{t}dB_{t}\right)
\left( \int_{n}^{n+1}\mathrm{sign}B_{t}dB_{t}\right) \right] .
\end{equation*}%
Using the formula for the expectation of the product of two divergence
integrals we obtain%
\begin{eqnarray*}
r(n) &=&\alpha _{H}\int_{a}^{a+1}\int_{n}^{n+1}\mathbb{E}\left( \mathrm{sign}%
B_{s}\mathrm{sign}B_{t}\right) \ |s-t|^{2H-2}dsdt \\
&&+4\alpha _{H}^{2}\int_{a}^{a+1}\int_{n}^{n+1}\int_{0}^{s}\int_{0}^{t}%
\mathbb{E}\left( \delta _{0}(B_{s})\delta _{0}(B_{t}\right) ) \\
&&\times \ |s-\sigma |^{2H-2}|\theta -t|^{2H-2}| d\sigma d\theta dsdt,
\end{eqnarray*}%
where $\alpha _{H}=H(2H-1)$. This formula can be proved by approximating the
function $\mathrm{sign}(x)$ by smooth functions and then taking the limit in
$L^{2}(\Omega )$. We have%
\begin{equation*}
\int_{0}^{t}|s-\sigma |^{2H-2}d\sigma =\frac{1}{2H-1}(s^{2H-1}+|s-t|^{2H-1}%
\mathrm{sign}(t-s)).
\end{equation*}%
Hence,
\begin{eqnarray*}
r(n) &=&\alpha _{H}\int_{a}^{a+1}\int_{n}^{n+1}\mathbb{E}\left( \mathrm{sign}%
B_{s}\mathrm{sign}B_{t}\right) \ |s-t|^{2H-2}dsdt \\
&&+4H^{2}\int_{a}^{a+1}\int_{n}^{n+1}\mathbb{E}\left( \delta
_{0}(B_{s})\delta _{0}(B_{t}\right) ) \\
&&\times (s^{2H-1}+(t-s)^{2H-1})(t^{2H-1}-(t-s)^{2H-1})dsdt \\
&=&a_{n}+b_{n}.
\end{eqnarray*}%
For the second term we have%
\begin{equation*}
b_{n}=\frac{4H^{2}}{2\pi }\int_{a}^{a+1}\int_{n}^{n+1}\frac{%
(s^{2H-1}+(t-s)^{2H-1})(t^{2H-1}-(t-s)^{2H-1})}{\left[ (st)^{2H}-\frac{1}{4}%
\left( t^{2H}+s^{2H}-|t-s|^{2H}\right) ^{2}\right] ^{1/2}}dsdt
\end{equation*}%
Therefore,
\begin{eqnarray*}
b_{n} &\leq &\frac{4H^{2}(2H-1)}{2\pi }\frac{\left(
(a+1)^{2H-1}+(n+1)^{2H-1}\right) (n-a-1)^{2H-2}}{\left[ (an)^{2H}-\frac{1}{4}%
\left( (n+1)^{2H}+(a+1)^{2H}-|n-a|^{2H}\right) ^{2}\right] ^{1/2}} \\
&\leq &Cn^{3H-3}.
\end{eqnarray*}
To estimate $a_{n}$, we have from (\ref{e.1.1})
\begin{eqnarray*}
{{\mathbb{E}}}\left[ {\mathrm{sign}}(B_{u}){\mathrm{sign}}(B_{v})\right]
&=&\sum_{k=0}^{\infty }\frac{4(2k)!}{(2k+1)^2 2\pi (k!2^{k})^{2}(uv)^{(2k+1)H}}%
(u^{2H}+v^{2H}-|u-v|^{2H})^{2k+1} \\
&\leq &C\frac{u^{2H}+v^{2H}-|u-v|^{2H}}{(uv)^{H}}\sum_{k=0}^{\infty }k^{-3/2}%
\frac{(u^{2H}+v^{2H}-|u-v|^{2H})^{2k}}{(uv)^{2kH}} \\
&\leq &C_{1}\frac{u^{2H}+v^{2H}-|u-v|^{2H}}{(uv)^{H}}.
\end{eqnarray*}%
Therefore
\begin{eqnarray*}
a_{n} &\leq &C_{2}\int_{a}^{a+1}\int_{n}^{n+1}|u-v|^{2H-2}\frac{%
u^{2H}+v^{2H}-|u-v|^{2H}}{(uv)^{H}}dvdu \\
&\leq &C_{3}n^{3H-3}\,.
\end{eqnarray*}%
As a consequence, if $H<2/3$, then$\ \sum_{n\geq 1}r(n)<\infty .$
\end{proof}

\section{General Dimension}

In this section we consider a $d$-dimensional fractional Brownian motion $%
B=\{(B_{t}^{1},\ldots ,B_{t}^{d}),t\geq 0\}$, with Hurst parameter $H>\frac{1%
}{2}$. Let $R_{t}=\ |B_{t}|$ be the fractional Bessel process associated to
the $d$-dimensional fBm $B$.

Suppose first that $H>\frac{1}{2}$. Fix a time interval $[0,T]$, and define
the derivative and divergence operators, $D^{\left( i\right) }$ and $\delta
^{\left( i\right) }$, with respect to each component $B^{\left( i\right) }$,
as in Section\ 2. We assume that the Sobolev spaces $\mathbb{D}_{i}^{1,p}$
include functionals of all the components of $B$ and not only of component $i
$. For each $p>1$, let $\mathbb{L}_{H,i}^{1,p}$ be the set of processes $%
u\in \mathbb{D}_{i}^{1,p}(\mathcal{H)}$ such that%
\begin{equation*}
\mathbb{E}\left[ \left\| u\right\| _{L^{1/H}\left( \left[ 0,T\right] \right)
}^{p}\right] +\mathbb{E}\left[ \left\| D^{(i)}u\right\| _{L^{1/H}\left( %
\left[ 0,T\right] ^{2}\right) }^{p}\right] <\infty \text{.}
\end{equation*}%
It has been proved in \cite{ChN} that $\ \left\{ \frac{B_{s}^{i}}{R_{s}}%
,s\in \left[ 0,T\right] \right\} $ $\ $belongs to the space $\mathbb{L}%
_{H,i}^{1,1/H}$ for each $i=1,\ldots ,d$ and
\begin{equation}
R_{t}=\sum_{i=1}^{d}\int_{0}^{t}\frac{B_{s}^{i}}{R_{s}}dB_{s}^{i}+H(d-1)%
\int_{0}^{t}\frac{s^{2H-1}}{R_{s}}ds.  \label{e10}
\end{equation}

In the case $H<\frac{1}{2}$ the following result holds.

\begin{proposition}
If  $H<\frac{1}{2}$, the process $\frac{B_{s}^{i}}{R_{s}}$ belongs to the
extended domain of the divergence operator $\mathrm{Dom}_{t}^{\ast }\delta
^{i}$ on any time interval $[0,t]$, and (\ref{e10}) holds.
\end{proposition}

\begin{proof}
For any test random variable $F\in \mathcal{S}_{\mathcal{H}}$ we have
\begin{eqnarray*}
&&\int_{0}^{t}\mathbb{E}(\frac{B_{s}^{i}}{R_{s}}K_{H}^{\ast ,a}K_{H}^{\ast
}D_{s}^{(i)}F)ds \\
&=&\lim_{\varepsilon \downarrow 0}\int_{0}^{t}\mathbb{E}(h_{\varepsilon
}^{\prime }(R_{s}^{2})B_{s}^{i})K_{H}^{\ast ,a}K_{H}^{\ast }D_{s}^{(i)}F)ds
\\
&=&\lim_{\varepsilon \downarrow 0}\ \mathbb{E}\left(
F\int_{0}^{t}h_{\varepsilon }^{\prime }(R_{s}^{2})B_{s}^{i}dB_{s}^{i}\right)
,
\end{eqnarray*}%
where, for any $\varepsilon >0$,
\begin{equation*}
h_{\varepsilon }(x)=\left\{
\begin{array}{cc}
\frac{3}{8}\sqrt{\varepsilon }+\frac{3}{4\sqrt{\varepsilon }}x-\frac{1}{%
8\varepsilon \sqrt{\varepsilon }}x^{2} & \text{ if }x<\varepsilon  \\
\sqrt{x} & \text{ if }x\geq \varepsilon
\end{array}%
\right. .
\end{equation*}%
We have $h_{\varepsilon }(x)\in C^{2}\left( \mathbb{R}\right) $ and $%
\underset{\varepsilon \rightarrow 0}{\lim }h_{\varepsilon }(x)=\sqrt{x}$ for
all $x\geq 0$. By It\^{o}'s formula for the fractional Brownian motion in
the case $H<\frac{1}{2}$, we have%
\begin{equation*}
h_{\varepsilon }(R_{t}^{2})-h_{\varepsilon
}(0)=\sum_{i=1}^{d}\int_{0}^{t}h_{\varepsilon }^{\prime
}(R_{s}^{2})B_{s}^{i}dB_{s}^{i}+J_{\varepsilon },
\end{equation*}%
where%
\begin{eqnarray*}
J_{\varepsilon } &=&H(d-1)\int_{0}^{t}\mathbf{1}_{\left\{ R_{s}^{2}\geq
\varepsilon \right\} }\frac{s^{2H-1}}{R_{s}}ds \\
&&+H\int_{0}^{t}\mathbf{1}_{\left\{ R_{s}^{2}<\varepsilon \right\} }\frac{1}{%
2\sqrt{\varepsilon }}\left[ 3d-(d+2)\frac{R_{s}^{2}}{\varepsilon }\right]
s^{2H-1}ds.
\end{eqnarray*}
\end{proof}

The process\
\begin{equation}
X_{t}=\sum_{i=1}^{d}\int_{0}^{t}\frac{B_{s}^{i}}{R_{s}}\delta B_{s}^{i}
\label{e2}
\end{equation}%
is $H$-self-similar, H\"{o}lder continuous of order $\alpha <H$, and it has
the same $\frac{1}{H}$-variation as the fractional Brownian motion.
Nevertheless, as we will show in the next proposition it is not a fractional
Brownian motion with Hurst parameter $H$.

For $h\in \mathcal{H}^{\otimes n}$, we denote
\begin{equation}
I_{j_{1},\ldots ,j_{n}}(h)=\int_{0<s_{1},\ldots ,s_{n}<t}h(s_{1},\ldots
,s_{n})dB_{s_{1}}^{j_{1}}\cdots dB_{s_{n}}^{j_{n}}\,,
\end{equation}

First we find the chaos expansion of $\sum_{i=1}^{d}%
\int_{0}^{t}f_{i}(B_{s})dB_{s}^{i}$, where $f_i$ $:\mathbb{R}^{d}\rightarrow \mathbb{R}$
are  smooth functions  with polynomial growth.

\begin{proposition}
The following chaos expansion holds for $Z_t=\sum_{i=1}^{d}%
\int_{0}^{t}f_{i}(B_{s})dB_{s}^{i}$
\begin{equation*}
Z_t=\sum_{i=1}^{d}\sum_{n=1}^{\infty }\sum_{1\leq j_{1},\ldots ,j_{n}\leq
d}I_{j_{1},\ldots ,j_{n}, i}\left( g_{j_{1},\ldots ,j_{n}}^{i}(s_{1},\ldots
,s_{n+1}\right)\,,
\end{equation*}%
where
\begin{eqnarray*}
g_{j_{1},\ldots ,j_{n}}^{i}(s_{1},\ldots ,s_{n+1}) &=&\frac{%
(-1)^{n}(s_{1}\vee \cdots \vee s_{n+1})^{-nH}}{(2\pi )^{d/2}} \\
&&\times \int_{{{\mathbb{R}}}^{d}}\left[ \frac{\partial ^{n}}{\partial
y_{j_{1}}\cdots \partial y_{j_{n}}}e^{-\frac{|y|^{2}}{2}}\right]
f_{i}(y(s_{1}\vee \cdots \vee s_{n+1})^{H})dy.
\end{eqnarray*}
\end{proposition}

\begin{proof}
Using Stroock's formula yields for each $i=1,\ldots ,d$%
\begin{equation*}
f_{i}(B_{s})=\sum_{n=0}^{\infty }\sum_{1\leq j_{1},\ldots ,j_{n}\leq d}\frac{%
1}{n!}I_{j_{1},\ldots ,j_{n}}\left( f_{j_{1},\ldots ,j_{n}}^{i}(s)\mathbf{1}%
_{[0,s]}^{\otimes n}\right) ,
\end{equation*}%
where%
\begin{eqnarray*}
f_{j_{1},\ldots ,j_{n}}^{i}(s) &=&\mathbb{E}(D^{j_{1}}\cdots
D^{j_{n}}(f_{i}(B_{s}))) \\
&=&\mathbb{E}(\frac{\partial ^{n}f_{i}}{\partial z_{j_{1}}\cdots \partial
z_{j_{n}}}(B_{s})) \\
&=&\frac{1}{(2\pi s^{2H})^{d/2}}\int_{\mathbb{R}^{d}}\frac{\partial ^{n}f_{i}%
}{\partial z_{j_{1}}\cdots \partial z_{j_{n}}}(z)e^{-\frac{|z|^{2}}{2s^{2H}}%
}dz \\
&=&\frac{s^{-nH}}{(2\pi )^{d/2}}\int_{\mathbb{R}^{d}}\frac{\partial ^{n}f_{i}%
}{\partial y_{j_{1}}\cdots \partial y_{j_{n}}}(ys^{H})e^{-\frac{|y|^{2}}{2y}%
}dy \\
&=&\frac{(-1)^{n}s^{-nH}}{(2\pi )^{d/2}}\int_{\mathbb{R}^{d}}f_{i}(ys^{H})%
\frac{\partial ^{n}}{\partial y_{j_{1}}\cdots \partial y_{j_{n}}}e^{-\frac{%
|y|^{2}}{2y}}dy.
\end{eqnarray*}%
Finally, the result follows from%
\begin{eqnarray*}
Z_t &=&\sum_{i=1}^{d}\int_{0}^{t}f_{i}(B_{s})dB_{s}^{i} \\
&=&\sum_{i=1}^{d}\sum_{n=0}^{\infty }\sum_{1\leq j_{1},\ldots ,j_{n}\leq
d}I_{j_{1},\ldots ,j_{n},i}\left( \mathrm{symm}\left( f_{j_{1},\ldots
,j_{n}}^{i}(s)\mathbf{1}_{[0,s]}^{\otimes n}\mathbf{1}_{[0,t]}(s)\right)
\right) \\
&=&\sum_{i=1}^{d}\sum_{n=0}^{\infty }\sum_{1\leq j_{1},\ldots ,j_{n}\leq
d}I_{j_{1},\ldots ,j_{n},i}\left( f_{j_{1},\ldots ,j_{n}}^{i}(s_{1}\vee
\cdots \vee s_{n+1})\prod_{i=1}^{n+1}\mathbf{1}_{[0,t]}^{{}}(s_{i})\right) ,
\end{eqnarray*}%
which completes the proof of the proposition.
\end{proof}

Now let $f_{i}(x)=\frac{x_{i}}{\sqrt{x_{1}^{1}+\cdots +x_{d}^{2}}}$. Then it
is easy to check $f_{i}(tx)=f_{i}(x)$ for all $t>0$. Hence, for such $f_{i}$%
, we have
\begin{equation*}
g_{j_{1},\ldots ,j_{n}}^{i}(s_{1},\ldots ,s_{n+1})=b_{j_{1},\cdots
,j_{n}}(s_{1}\vee \cdots \vee s_{n+1})^{-nH},
\end{equation*}
where%
\begin{equation*}
b_{i,j_{1},\ldots ,j_{n}}=\frac{(-1)^{n}}{(2\pi )^{d/2}}\int_{{{\mathbb{R}}}%
^{d}}\left[ \frac{\partial ^{n}}{\partial y_{j_{1}}\cdots \partial y_{j_{n}}}%
e^{-\frac{|y|^{2}}{2}}\right] f_{i}(y)dy.
\end{equation*}%
Then the chaos expansion of $f_{i}(B_{t})$ is given by
\begin{equation*}
f_{i}(B_{t})=\sum_{n=0}^{\infty }\sum_{1\leq j_{1},\ldots ,j_{n}\leq
d}b_{i,j_{1},\ldots ,j_{n}}t^{-nH}\int_{\left\{ 0<s_{1}<\cdots
<s_{n}<t\right\} }dB_{s_{1}}^{j_{1}}\cdots dB_{s_{n}}^{j_{n}}\,,
\end{equation*}%
and the chaos expansion of the divergence of this process is
\begin{eqnarray*}
\int_{0}^{t}f_{i}(B_{s})dB_{s}^{i} &=&\sum_{n=1}^{\infty }\sum_{1\leq
j_{1},\ldots ,j_{n}\leq d}b_{i,j_{1},\ldots ,j_{n}} \\
&&\times \int_{\left\{ 0<s_{1},\ldots ,s_{n+1}<t\right\} }(s_{1}\vee \cdots
\vee s_{n+1})^{-nH}dB_{s_{1}}^{j_{1}}\cdots
dB_{s_{n}}^{j_{n}}dB_{s_{n+1}}^{i}.
\end{eqnarray*}

Using these results we can prove the following proposition.

\begin{proposition}
The process $Y=\{Y_{t},t\geq 0\}$ defined in (\ref{e2}) is not a fractional
Brownian motion.
\end{proposition}

\begin{proof}
If $Y_{t}=\sum_{i=1}^{d}\int_{0}^{t}f_{i}(B_{s})dB_{s}^{i}$ is a fractional
Brownian motion, then as in previous section one can show that
\begin{equation*}
Z_{t}=\sum_{i=1}^{d}\int_{0}^{t}\eta _{H}(t,s)f_{i}(B_{s})dB_{s}^{i}
\end{equation*}%
is a classical Brownian motion. But
\begin{equation*}
Z_{t}=\sum_{i=1}^{d}\sum_{n=1}^{\infty }\sum_{1\leq j_{1},\ldots ,j_{n}\leq
d}b_{i,j_{1},\ldots ,j_{n}}\int_{\left\{ 0<s_{1},\ldots ,s_{n+1}<t\right\}
}h(t,s_{1},\ldots ,s_{n+1})dB_{s_{1}}^{j_{1}}\cdots
dB_{s_{n}}^{j_{n}}dB_{s_{n+1}}^{i}\,,
\end{equation*}%
where
\begin{equation*}
h(t,s_{1},\ldots ,s_{n+1})=\eta (t,s_{1}\vee s_{2}\vee \cdots \vee
s_{n+1})(s_{1}\vee s_{2}\vee \cdots \vee s_{n+1})^{-nH}\,.
\end{equation*}%
In a similar way to one dimensional case, one can show that 
$\{ Z_t, t\ge 0\}$   is not a
martingale
\end{proof}

\begin{proposition}
Let $Y=\{Y_{t},t\geq 0\}$ be the process defined in (\ref{e2}) If $H\geq 2/3$%
, then $Y_{t}$ is long range dependent and if $1/2<H<2/3$, then $Y_{t}$ is
not long range dependent.
\end{proposition}

\begin{proof}
Set%
\begin{equation*}
\rho _{n}=\mathbb{E}\left[ \left( \sum_{i=1}^{d}\int_{0}^{1}\frac{B_{s}^{i}}{%
|B_{s}|}\delta B_{s}^{i}\right) \left( \sum_{i=1}^{d}\int_{n}^{n+1}\frac{%
B_{s}^{i}}{|B_{s}|}\delta B_{s}^{i}\right) \right] .
\end{equation*}%
By the formula for the expectation of the product of two divergence
integrals we can write%
\begin{eqnarray*}
\rho _{n} &=&\sum_{i,j=1}^{d}\alpha _{H}\int_{0}^{1}\int_{n}^{n+1}\mathbb{E}%
\left( \frac{B_{s}^{i}B_{t}^{i}}{|B_{s}||B_{t}|}\right) \ |t-s|^{2H-2}dsdt \\
&&+\sum_{i,j=1}^{d}\alpha
_{H}^{2}\int_{0}^{1}\int_{n}^{n+1}\int_{0}^{t}\int_{0}^{s}\mathbb{E}\left(
D_{\theta }^{j}\left( \frac{B_{s}^{i}}{|B_{s}|}\right) D_{\sigma }^{i}\left(
\frac{B_{t}^{j}}{|B_{t}|}\right) \right) \\
&&\times |\theta -t|^{2H-2}|\sigma -s|^{2H-2}d\theta d\sigma dsdt \\
&:&=\rho _{n}^{1}+\rho _{n}^{2}.
\end{eqnarray*}%
In order to estimate the term $\rho _{n}^{1}$ we make use of the orthogonal
decomposition
\begin{equation*}
B_{t}=\frac{R(t,s)}{s^{2H}}B_{s}+\beta _{s,t}Y,
\end{equation*}%
where $\ $
\begin{equation*}
\beta _{s,t}^{2}=\frac{(st)^{2H}-R(t,s)^{2}}{s^{2H}},
\end{equation*}%
and $Y$ is a $d$-dimensional standard normal random variable independent of $%
B_{s}$ . \ Set%
\begin{equation*}
\lambda _{st}=\frac{R(t,s)}{\beta _{s,t}s^{2H}}=\frac{\frac{1}{2}\left(
t^{2H}+s^{2H}-|t-s|^{2H}\right) }{s^{H}\left[ (st)^{2H}-\frac{1}{4}\left(
t^{2H}+s^{2H}-|t-s|^{2H}\right) ^{2}\right] ^{1/2}}
\end{equation*}%
As $t$ tends to infinity and $s$ belongs to $(0,1)$, the term $\lambda _{st}$
behaves as $Hs^{-2H}t^{H-1}$. Hence, by Lemma \ref{lem1}%
\begin{eqnarray*}
\mathbb{E}\left( \frac{\left\langle B_{s},B_{t}\right\rangle }{|B_{s}||B_{t}|%
}\right) &=&\mathbb{E}\left( \frac{\left\langle B_{s},\lambda
_{st}B_{s}+Y\right\rangle }{|B_{s}||\lambda _{st}B_{s}+Y|}\right) \\
&=&\mathbb{E}\left( \frac{\left\langle B_{1},s^{H}\lambda
_{st}B_{1}+Y\right\rangle }{|B_{1}||s^{H}\lambda _{st}B_{1}+Y|}\right) \\
&=&s^{H}\lambda _{st}\mathbb{E}\left( \frac{\left| B_{1}\right| ^{2}\left|
Y\right| ^{2}-\left\langle B_{1},Y\right\rangle ^{2}}{\left| B_{1}\right|
\left| Y\right| ^{3}}\right) +o(t^{H-1}).
\end{eqnarray*}%
Hence,$\mathbb{\ }$%
\begin{equation*}
\mathbb{E}\left( \frac{\left\langle B_{s},B_{t}\right\rangle }{|B_{s}||B_{t}|%
}\right) \thickapprox HCs^{-H}t^{H-1},
\end{equation*}%
where
\begin{equation*}
C=\mathbb{E}\left( \frac{\left| B_{1}\right| ^{2}\left| Y\right|
^{2}-\left\langle B_{1},Y\right\rangle ^{2}}{\left| B_{1}\right| \left|
Y\right| ^{3}}\right) >0.
\end{equation*}%
This implies that the term $\rho _{n}^{1}$ behaves as $n^{3H-3}$.

For the term $\rho _{n}^{2}$ we have%
\begin{eqnarray*}
\rho _{n}^{2} &=&H^{2}\sum_{i,j=1}^{d}\int_{0}^{1}\ \int_{n}^{n+1}\mathbb{E}%
\left( \ \left( \frac{\delta _{ij}}{|B_{s}|}-\frac{B_{s}^{i}B_{s}^{j}}{%
|B_{s}|^{3}}\right) \ \left( \frac{\delta _{ij}}{|B_{t}|}-\frac{%
B_{t}^{i}B_{t}^{j}}{|B_{t}|^{3}}\right) \right) \\
&&\times (s^{2H-1}+(t-s)^{2H-1})(t^{2H-1}-(t-s)^{2H-1})dsdt \\
&=&H^{2}\int_{0}^{1}\ \int_{n}^{n+1}\mathbb{E}\left( \ \frac{d}{%
|B_{s}||B_{t}|}-\frac{|B_{t}|^{2}}{|B_{s}||B_{t}|^{3}}-\frac{|B_{s}|^{2}}{%
|B_{s}|^{3}|B_{t}|}+\frac{\left\langle B_{s},B_{t}\right\rangle ^{2}}{%
|B_{s}|^{3}|B_{t}|^{3}}\right) \\
&&\times (s^{2H-1}+(t-s)^{2H-1})(t^{2H-1}-(t-s)^{2H-1})dsdt \\
&=&H^{2}\int_{0}^{1}\ \int_{n}^{n+1}\mathbb{E}\left( \ \frac{d-2}{%
|B_{s}||B_{t}|}+\frac{\left\langle B_{s},B_{t}\right\rangle ^{2}}{%
|B_{s}|^{3}|B_{t}|^{3}}\right) \\
&&\times (s^{2H-1}+(t-s)^{2H-1})(t^{2H-1}-(t-s)^{2H-1})dsdt.
\end{eqnarray*}%
The term%
\begin{equation*}
\mathbb{E}\left( \ \frac{1}{|B_{s}||B_{t}|}\left( d-2+\frac{\left\langle
B_{s},B_{t}\right\rangle ^{2}}{|B_{s}|^{2}|B_{t}|^{2}}\right) \right)
\end{equation*}%
behaves as $Kt^{-H}$ as $t$ tends to infinity, where%
\begin{equation*}
K=\mathbb{E}\left( \ \frac{1}{|B_{1}||Y|}\left( d-2+\frac{\left\langle
B_{1},Y\right\rangle ^{2}}{|B_{1}|^{2}|Y|^{2}}\right) \right) >0.
\end{equation*}%
Hence, the term $\rho _{n}^{2}$ behaves also as $n^{3H-3}$. This completes
the proof taking into account that the constants $C$ and $K$ are positive.
\end{proof}

\begin{lemma}
\label{lem1} Let $X$ and $Y$ be independent $d$-dimensional standard normal
random variables. Then as $\varepsilon $ tends to zero we have
\begin{equation*}
\mathbb{E}\left( \frac{\left\langle X,\varepsilon X+Y\right\rangle }{\left|
X\right| \left| \varepsilon X+Y\right| }\right) =\varepsilon \mathbb{E}%
\left( \frac{\left| X\right| ^{2}\left| Y\right| ^{2}-\left\langle
X,Y\right\rangle ^{2}}{\left| X\right| \left| Y\right| ^{3}}\right)
+o(\varepsilon ).
\end{equation*}
\end{lemma}

\begin{proof}
We have%
\begin{eqnarray*}
\mathbb{E}\left( \frac{\left\langle X,\varepsilon X+Y\right\rangle }{\left|
X\right| \left| \varepsilon X+Y\right| }\right) &=&\mathbb{E}\left( \frac{%
\left\langle X,\varepsilon X+Y\right\rangle }{\left| X\right| \left|
\varepsilon X+Y\right| }-\frac{\left\langle X,Y\right\rangle }{\left|
X\right| \left| Y\right| }\right) \\
&=&\mathbb{E}\left( \frac{\left\langle X,\varepsilon X+Y\right\rangle \left|
Y\right| -\left\langle X,Y\right\rangle \left| \varepsilon X+Y\right|
\mathbb{\ }}{\left| X\right| \left| \varepsilon X+Y\right| \left| Y\right| }%
\ \right) \\
&=&\mathbb{E}\left( \frac{\varepsilon \left| X\right| ^{2}\left| Y\right|
-\varepsilon \left\langle X,Y\right\rangle ^{2}/\left| Y\right| \mathbb{\ +}o%
\mathbb{(\varepsilon )}}{\left| X\right| \left| \varepsilon X+Y\right|
\left| Y\right| }\ \right) ,
\end{eqnarray*}%
and that yields the desired estimation.
\end{proof}

\end{document}